\documentclass[10pt,notitlepage,twoside,a4paper]{amsart}
 \usepackage{mathrsfs}

\usepackage{amsmath,amssymb,enumerate}

\usepackage{epsfig,fancyhdr,color}

\usepackage{amssymb}
\usepackage{amsmath,amsthm}
\usepackage{latexsym}
\usepackage{amscd}
\usepackage{psfrag}
\usepackage{graphicx}
\usepackage[latin1]{inputenc}
\usepackage[all]{xy}

\newtheorem{theorem}{\rm\bf Theorem}[section]
\newtheorem{proposition}[theorem]{\rm\bf Proposition}
\newtheorem{lemma}[theorem]{\rm\bf Lemma}

\newtheorem*{theorem 1}{\rm\bf Proposition 1}
\newtheorem*{theorem 2}{\rm\bf Proposition 2}

\theoremstyle{definition}

\theoremstyle{remark}
\newtheorem{remark}[theorem]{\rm\bf Remark}

\def\interieur#1{\mathord{\mathop{\kern 0pt #1}\limits^\circ}}

\title[ Quasiisometry]{Almost-isometry between Teichm\"{u}ller metric and
length-spectra metric on moduli space}

\author{Lixin Liu}
\address{L. Liu, Department of Mathematics, Sun Yat-sen (Zhongshan) University, 510275, Guangzhou, P. R. China}
\email{mcsllx@mail.sysu.edu.cn}

\author{Weixu Su}
\address{W. Su, Department of Mathematics, Sun Yat-sen (Zhongshan) University, 510275, Guangzhou, P. R. China}
\email{su023411040@163.com}

\thanks{The work is partially supported by NSFC : 10871211}

\date{\today}

\begin{document}

\begin{abstract}
We prove an analogue of Farb-Masur's theorem that the length-spectra
metric on moduli space is  ``almost isometric'' to a simple model $\mathcal {V}(S)$
which is induced by the cone metric over the complex of curves. As
an application, we know that the Teichm\"{u}ller metric and the
length-spectra metric are ``almost isometric'' on moduli space,
while they are not even quasi-isometric on Teichm\"{u}ller space.
\bigskip



\end{abstract}
\maketitle
\section{Introduction}\label{intro}

         Let $S = S_{g,n}$ be an oriented  surface
of genus g with n punctures. We assume that $3g-3+n \geq 1 $. Let
$\mathcal{T}(S)$ denote the Teichm\"{u}ller space of equivalence
classes of marked Riemann surfaces on $S$.  A marked Riemann surface
is a pair $(X,f)$ where $X$ is a Riemann surface, considered as a surface with
either a conformal structure or hyperbolic metric on $S$, and $f:S\rightarrow
X$ is an orientation preserving homeomorphism. Two markings $(X_i,f_i),i=1,2$ are equivalent
if and only if there exists a conformal map $h:X_1\rightarrow X_2$ which is
homotopic to $f_2 \circ f_1^{-1}$. In the following, we will often
denote $(X,f) \in \mathcal {T}(S)$ by $X$, without explicit
reference to the marking or the equivalence relation.

There are several natural metrics on $\mathcal{T}(S)$. In this paper
we will consider the Teichm\"{u}ller metric $d_{\mathrm{Teich}}$ and the
length-spectra metric $d_{\mathrm{ls}}$. Both of the two metrics are complete Finsler metrics.
For the Teichm\"uller metric, this
is a classical well-known result.
The fact that the length spectrum metric is Finsler follows from the fact that this
metric is a symmetrization of Thurston's asymmetric metric, which is Finsler, which was proved by Thurston in \cite{Thurston}.
In this paper, the Finsler property of the length-spectra metric to used to show that any two points 
in the Teichm\"uller metric can be connected by a geodesic. 

Recall that for marked conformal structures $X_1, X_2\in
\mathcal{T}(S)$, the Teichm\"{u}ller metric is defined by
$$d_{\mathrm{Teich}}(X_1,X_2) = \frac{1}{2}\log K$$
where $K\geq 1$ is the least number such that there exists a
K-quasiconformal map between the marked structures $X_1$ and $X_2$.
Teichm\"{u}ller's Theorem states that there exists a unique extremal
quasiconformal map realizing $d_{\mathrm{Teich}}(X_1,X_2)$. See Abikoff
\cite{Abikoff} for details.

The length-spectra metric (also called the Lipschitz metric \cite{CR} )
on $\mathcal {T}(S)$ is defined by
$$d_{\mathrm{ls}}(X_1,X_2) = \frac{1}{2}\log \max\{K_1, K_2\}$$
where $K_1\geq 1$ is the least number such that there is a global
$K_1$-Lipschitz homeomorphism between the marked hyperbolic
metrics $X_1$ and $X_2$,
  and  where $K_2\geq 1$ is the least number
such that there is a global $K_2$-Lipschitz  homeomorphism between
the marked hyperbolic metrics $X_2$ and $X_1$. Note that there are
two asymmetric metrics (called Thurston's asymmetric metrics) defined by $d_1(X_1,X_2)=\frac{1}{2}\log K_1$
and $d_2(X_1,X_2)= \frac{1}{2}\log K_2$. Thurston \cite{Thurston}
showed that the extremal Lipschitz maps exist.

The mapping class group $\mathrm{Mod}(S)$ is the group of homotopy classes of
orientation-preserving homeomorphisms of $S$. This group acts
properly discontinuously and isometrically on $(\mathcal{T}(S),
d_{\mathrm{Teich}})$ and $(\mathcal{T}(S), d_{\mathrm{ls}})$, thus inducing two
metrics $d_\mathcal {T}$ and $d_\mathcal {L}$ on the quotient moduli
space $\mathcal{M}(S):= \mathcal {T}(S)/\mathrm{Mod}(S)$. Let $\pi: \mathcal{T}(S) \to \mathcal{M}(S)$
be the natural projection.

Kerckhoff \cite{Kerckhoff} discovered an elegant and useful formula
to compute the Teichm\"{u}ller metric in terms of extremal length:
$$d_{\mathrm{Teich}}(X_1,X_2) = \frac{1}{2}\log \sup_\gamma \frac{\mathrm{Ext}_{X_1}(\gamma)}{\mathrm{Ext}_{X_2}(\gamma)}$$
where the supremum is taken over all isotopy classes of essential
(neither homotopic to a point nor to a puncture) simple closed
curves on $S$.  The  extremal length of $\gamma$ in $X$, denoted by
 $\mathrm{Ext}_X(\gamma)$, is
 defined by
 $$\mathrm{Ext}_X(\gamma) := \sup_\rho \frac{{L_{\rho}(\gamma)}^2}{\mathrm{Area}_\rho}$$
 where the supremum is taken over all conformal metrics $\rho$ on $X$ of finite
 positive area.

On the other hand, it was shown by Thurston \cite{Thurston} that the minimal
Lipschitz constant is given by the ratios of hyperbolic length:
$$K_1 = \sup_\gamma\frac{l_{X_2}(\gamma)}{l_{X_1}(\gamma)}, K_2 = \sup_\gamma\frac{l_{X_1}(\gamma)}{l_{X_2}(\gamma)}$$
where the supremum is taken over all isotopy classes of essential
simple closed curves on $S$. Thus the length-spectra
metric is given by
\begin{equation}\label{equ:ls}
d_{\mathrm{ls}}(X_1,X_2) = \max \{\frac{1}{2}\log
\sup_\gamma\frac{l_{X_2}(\gamma)}{l_{X_1}(\gamma)},\frac{1}{2}\log
\sup_\gamma\frac{l_{X_1}(\gamma)}{l_{X_2}(\gamma)}\}.
\end{equation}

It is of interest to study the relation between the metrics
$d_{\mathrm{Teich}}$ and $d_{\mathrm{ls}}$. The following lemma of Wolpert
\cite{Wolpert} implies that
$d_{\mathrm{ls}}(X_1,X_2)\leq d_{\mathrm{Teich}}(X_1,X_2)$.
\begin{lemma} \label{lem:four}
For any K-quasiconformal mapping $f$ from $X_1$ to $X_2$ and any
simple closed curves $\gamma$, we have
$$\frac{l_{X_2}(f(\gamma))}{l_{X_1}(\gamma)}\leq K.$$
 \end{lemma}
It was shown by Li \cite{L1} that $d_{\mathrm{Teich}}$ and $d_{\mathrm{ls}}$ induce
the same topology on $\mathcal{T}(S)$.  Moreover,
$d_{\mathrm{ls}}(X_1,X_2)\leq d_{\mathrm{Teich}}(X_1,X_2)\leq 2d_{\mathrm{ls}}(X_1,X_2) +
C(X_1)$, where $C(X_1)$ is a constant depending on $X_1$. The proof
can be shown by considering the ratios of extremal length and square of
hyperbolic length $\frac{\mathrm{Ext}_X(\gamma)}{l_X^2(\gamma)}$, which
defines a function on the compact space of projective measured
foliations.

However, Li \cite{L2} also  proved that $d_{\mathrm{Teich}}$ and $d_{\mathrm{ls}}$ are
not metrically equivalent; that is, there is no constant $C>0$ such
that $d_{\mathrm{ls}}(X_1,X_2)\leq d_{\mathrm{Teich}}(X_1,X_2)\leq Cd_{\mathrm{ls}}(X_1,X_2)$
for any $X_1$ and $X_2$ in $\mathcal {T}(S)$. In particular, Choi
and Rafi \cite{CR} (see also Liu, Sun and Wei \cite{LSW}) showed
that although the two metrics are quasi-isometric to each other on
the thick part of $\mathcal{T}(S)$, there are sequences
$X_n,Y_n,n=1,2,\cdots$ in the thin part of $\mathcal{T}(S)$, such
that $d_{\mathrm{ls}}(X_n,Y_n)\rightarrow 0$, while
$d_{\mathrm{Teich}}(X_n,Y_n)\rightarrow \infty$.

In fact, all the well-known examples that illustrate the divergence
of $d_{\mathrm{Teich}}$ and $d_{\mathrm{ls}}$ are constructed by Dehn twists. Choi and
Rafi \cite{CR} also noticed the fact that for any two points $X_1,
X_2$ in the thin part of $\mathcal{T}(S)$, if they have no short
curves in common, then $d_{\mathrm{Teich}}(X_1,X_2)$ is comparable to
$d_{\mathrm{ls}}(X_1,X_2)$. These give evidences that $d_{\mathrm{Teich}}$ and
$d_{\mathrm{ls}}$ may be quasi-isometric on the moduli space
$\mathcal{M}(S)$.

For more recent progress on the length-spectra metric and
Teichm\"uller metric in Teichm\"uller space, please see \cite{LP},
\cite{LPST1} and \cite{LPST2}.

Recently, Farb and Masur \cite{FM} studied the large-scale geometry of
moduli space. They built an ``almost isometric'' simplicial model
for $\mathcal{M}(S)$ with the Teichm\"{u}ller metric, from which
they determine the tangent cone at infinity of $\mathcal{M}(S)$.
Their result can be seen as a step in providing a "reduction theory"
for the action of $\mathrm{Mod}(S)$ on $\mathcal {T}(S)$, in analogy with the
case of locally symmetric spaces. See Farb and Masur \cite{FM},
Leuzinger \cite{Leuzinger} for details.

Let $\mathcal{C}(S)$ be the complex of curves on $S$.  This was introduced by Harvey \cite{Harvey} as an analogue in the context of Teichm\"uller
space of the Tits building associated to an arithmetic group.The vertices of $\mathcal {C}(S)$ are the free
isotopy classes of essential simple closed curves
on $S$, and a k-simplex consist of $k+1$ isotopy classes of
mutually disjoint essential simple closed curves. Note that $\mathcal
{C}(S)$ is a simplicial complex of dimension $(d-1)$, where $d =
3g-3+n$. A $(d-1)$-simplex is called a maximal simplex. Every
simplex is the face of a maximal simplex. While $\mathcal{C}(S)$
is locally infinite, the mapping class group $\mathrm{Mod}(S)$ acts on
$\mathcal {C}(S)$ and the quotient $\mathcal {C}(S)/\mathrm{Mod}(S)$ is a
finite orbi-complex. See \cite{Harvey} for reference.

Denote $\widetilde{\mathcal{V}}(S)$ to be the topological cone
$$\widetilde{\mathcal{V}}(S):= \frac{[0,\infty)\times \mathcal {C}(S)}{\{0\}\times \mathcal {C}(S)}.$$
For each maximal simplex $\sigma$ of $\mathcal {C}(S)$, we will
think of the cone over $\sigma$ as an orthant with coordinates
$(x_1,\cdots,x_d)\in \mathbb{R}_{\geq 0}^d
$. This orthant is endowed with the standard sup
metric:
$$d((x_1,\cdots,x_d),(y_1,\cdots,y_d)):= \frac{1}{2}\max_{1\leq i\leq
d}|x_i-y_i|.$$

Since the metrics on the cones on any two such maximal simplices
agree on (the cone on) any common face, we can endow
$\widetilde{\mathcal{V}}(S)$ with the corresponding path metric. The
natural action of $\mathrm{Mod}(S)$ on $\mathcal {C}(S)$ induces an isometric
action on $\widetilde{\mathcal{V}}(S)$, thus induces a well-defined
metric $d_\mathcal{V}$ on the quotient
$$\mathcal{V}(S):= \widetilde{\mathcal{V}}(S)/\mathrm{Mod}(S).$$

Given $C\geq 0$ and $\lambda\geq 1$, a map $f: X\rightarrow Y$ is
called a $(\lambda,C)$ quasi-isometry
if
$$\frac{1}{\lambda}d_X(x,y) - C\leq d_Y(f(x),f(y))\leq \lambda d_X(x,y) + C$$
for any $x,y\in X$, and the $C$-neighborhood of $f(X)$ in $Y$ is
all of $Y$.  A $(1,C)$ quasi-isometry
is called an ``almost isometry''.

The following theorem is the main result of Farb and Masur
\cite{FM}, which provides a simple and geometric model for
$\mathcal{M}(S)$.
\begin{theorem}\label{th:FM}
There is a map $\Psi: (\mathcal{V}(S),d_{\mathcal{V}})
\to (\mathcal{M}(S),d_\mathcal {T})$ which is an almost
isometry. That is, there is a constant $D$ that depends on $S$ such that:
\begin{enumerate}[$\bullet$]
\item $|d_\mathcal{V}(x,y)-d_\mathcal {T}(\Psi(x),\Psi(y))|\leq D$
for each $x,y\in \mathcal{V}(S)$, and
\item the $D$-neighborhood of
$\Psi(\mathcal{V}(S))$ in $(\mathcal{M}(S),d_\mathcal {T})$ is all of
$\mathcal{M}(S)$.
\end{enumerate}
 \end{theorem}

As a corollary of Theorem \ref{th:FM}, the tangent cone at
infinity of $\mathcal{M}(S)$ with the Teichm\"{u}ller metric is
isometric to $\mathcal{V}(S)$ and has dimension $d$.

The main goal of our article is to show that:

\begin{theorem}\label{th:LL}
Endow $\mathcal{M}(S)$ with the length-spectrum metric $d_\mathcal
{L}$, then the map $\Psi: \mathcal{V}(S) \to \mathcal{M}(S)$
in Theorem $1.2$ is an almost isometry. That is, there is a constant
$D'$ that depends on $S$ such that:
\begin{enumerate}[$\bullet$]
\item $|d_\mathcal{V}(x,y)-d_\mathcal {L}(\Psi(x),\Psi(y))|\leq D'$
for each $x,y\in \mathcal{V}(S)$, and
\item the $D'$-neighborhood of
$\Psi(\mathcal{V}(S))$ in $\mathcal{M}(S)$ is all of
$\mathcal{M}(S)$.
\end{enumerate}
 \end{theorem}

As a result, the tangent cone at infinity of $\mathcal{M}(S)$ with
the length-spectrum metric is also isometric to $\mathcal{V}(S)$.
Since quasi-isometry is an equivalence relation between metric
spaces, it is clear that Theorem 1.2 and Theorem 1.3 together imply
that
\begin{theorem}\label{th:LL}
The Teichm\"{u}ller metric and the length-spectrum metric are almost
isometric on $\mathcal{M}(S)$. That is, there is a constant $D''$
that depends on $S$, such that
$$d_\mathcal {L}(X_1,X_2)\leq d_\mathcal {T}(X_1,X_2)\leq d_\mathcal {L}(X_1,X_2) + D''$$
 for any $X_1, X_2\in \mathcal{M}(S)$.
 \end{theorem}

We will give a proof of Theorem $1.3$ in Section $3$. The method of
the proof of the theorem relies on Minsky's product theorem
\cite{Minsky}. Another ingredient of the proof is that any two
points of $\Psi(\mathcal {V}(S))$ can be joined by a quasi-geodesic
that lies in $\Psi(\mathcal {V}(S))$ and enters each simplex of
$\Psi(\mathcal {V}(S))$ at most once, as
observed by Farb and Masur \cite{FM}.\\
\textbf{Acknowledgements.}  The authors would like to thank
Professor Benson Farb and Professor Howard Masur for their helpful
suggestions and their interest in this project. They also thank the referee for his (or her) suggestions.

 \section{The map $\Psi$ }\label{s2}
We will define the map $\Psi:\mathcal{V}(S) \rightarrow
\mathcal{M}(S)$ as constructed by Farb-Masur \cite{FM}.

Let us first fix some notations.
Given a maximal simplex $\sigma\in \mathcal{C}(S)$, the cone over $\sigma$ is denoted by
$\bar\Delta(\sigma)$ and the open cone over $\sigma$ (with no $x_i=0$) is denoted by $\Delta(\sigma)$.
Let $\mathrm{Mod}(\sigma)$ be the subgroup of $\mathrm{Mod}(S)$ that fixes $\sigma$. It acts on $\Delta(\sigma)$
with finite orbit. Let $\Lambda(\sigma)$ be a sector inside $\Delta(\sigma)$ which is a fundamental domain for the
action of $\mathrm{Mod}(\sigma)$.

Now fix a maximal simplex $\sigma$ of $\mathcal{C}(S)$. $\sigma$ is
represented by a maximal collection of disjoint simple closed curves
$\{\alpha_1,\cdots,\alpha_d\}$. The choice of curves determines a
set of Fenchel-Nielsen coordinates on $\mathcal{T}(S)$, where a
point $X\in \mathcal{T}(S)$ is given by coordinates as following:
$$X\rightarrow (l_X(\alpha_1),\cdots,l_X(\alpha_d),\theta_1(X),\cdots,\theta_d(X))$$
where $l_X(\alpha_i)$ is the length of $\alpha_i$ with respect to
the hyperbolic metric $X$, and  $\theta_i$ are the so-called "twist
coordinates".

We now define a map $\widetilde{\Psi}:
\widetilde{\mathcal{V}}(S)\rightarrow \mathcal{M}(S)$ in the following way.
First $\widetilde{\Psi}$ is restricted
to the sector $\Lambda(\sigma)$ to be:
\begin{equation}\label{equ:psi}
\widetilde{\Psi}(x_1,\cdots,x_d) = \pi (X)
\end{equation}
where $X$ is a point of $\mathcal{T}(S)$ with $\ell_X(\alpha_i)=\epsilon_0e^{-x_i}$ and with twist coordinates
all equal to $0$.
Here the constant $\epsilon_0 = \epsilon_0(S)$ is a sufficiently
small constant such that for any hyperbolic surface $X$ homeomorphic
to $S$, any two simple closed curves $\alpha , \beta$ with length
not larger than $\epsilon_0$ are disjoint. Note that
$\epsilon_0e^{-x_i}\leq \epsilon_0$ and the image of
$\widetilde{\Psi}$ lies on $\Omega_\sigma(\epsilon_0)$,
where
$$\Omega_\sigma(\epsilon_0)=\{X\in \mathcal{T}(S): \ell_X(\alpha_i)<\epsilon, \ \mathrm{for} \ \mathrm{each} \ i=1, \cdots, d\}.$$
We use the action of $\mathrm{Mod}(\sigma)$ to extend $\widetilde{\Psi}$ from $\Lambda(\sigma)$ to $\Delta(\sigma)$. Since there are
finitely many collection of maximal simplices that represent all combinatorial types, we define the map $\widetilde{\Psi}$ for each
maximal cone in the finite collection. Then we use the action of $\mathrm{Mod}(S)$ to extend $\widetilde{\Psi}$ to the open cones on all
maximal simplices by having it be constant on orbits.

To define the map $\widetilde{\Psi}$ on $\widetilde{\mathcal{V}}(S)$, we have to define it on the cone over any simplex $\tau$ of $\mathcal{C}(S)$
that is not maximal. Since a simplex is a face of some maximal simplex (maybe not unique), we can choose some maximal simplex $\sigma$ that containing
$\tau$. The cone over $\tau$ is given by the coordinates $(x_1, \cdots, x_d)$ for the cone over $\sigma$. We define the map $\widetilde{\Psi}$
on the cone over $\tau$ via the equation $(\ref{equ:psi})$ by setting the coordinates $x_i$ corresponding curves in $\sigma-\tau$ to be $0$.
It follows that $\widetilde{\Psi}$ is $\mathrm{Mod}(S)$-invariant and thus induces a map $\Psi: \mathcal{V}(S) \to \mathcal{M}(S)$.

It is noticed that $\Psi$ is not continuous in general because of the choices made at a face of a maximal simplex.
Nevetheless we know that the jump in the function at any face is uniformly bounded.
Such an argument for the Teichm\"uller metric is proved by Farb and Masur \cite{FM}.
Moreover, they proved that the map $\Psi$ satisfies the
condition of Theorem 1.2, that is, $\Psi
:(\mathcal{V}(S),d_{\mathcal{V}})\rightarrow
(\mathcal{M}(S),d_\mathcal {T})$ is almost onto and almost
isometric.
To prove an analogue for $(\mathcal{M}(S),d_\mathcal
{L})$, we will show that the propositions that were used in the proof of
Theorem \ref{th:FM} are also applied to the length-spectra metric.

Let $\sigma=\{\alpha_1,\cdots,\alpha_d\}$ be a maximal simplex.
Following Minsky \cite{Minsky}, we change the Fenchel-Nielsen coordinates to
$$X\rightarrow
(\theta_1(X),\frac{1}{l_X(\alpha_1)},\cdots,\theta_d(X),\frac{1}{l_X(\alpha_d)})\in (\textbf{H}^2)^d.$$
We give $\textbf{H}^2$ the Poincar\'{e} metric $ds^2 = \frac{dx^2 +
dy^2}{4y^2}$ and endow $(\textbf{H}^2)^d$ with the sup metric.

We need the following lemma, which is a special case of the product
region theorem of Minsky \cite{Minsky}.

 \begin{lemma} \label{le:M}
With the notations above, there exists a constant $C$ depending
on $\epsilon_0$, such that for any $X, Y \in \Omega_\sigma(\epsilon_0)$.
$$|d_{\mathrm{Teich}}(X,Y) - \sup_{i = 1, \cdots, d} \{d_{\textbf{H}^2}((\theta_i(X),\frac{1}{l_X(\alpha_i)}),(\theta_i(X),\frac{1}{l_X(\alpha_i)}))\}| \leq C.$$
\end{lemma}

The following proposition is due to Farb and Masur \cite{FM},
which shows that $\Psi$ is almost onto. We include a proof here
for completeness.
 \begin{proposition}
There is a constant $C_1 = C_1(S)$, such that for any $X\in
\mathcal{M}(S)$, there exist a $Z\in {\Psi(\mathcal{V}(S))}$ such
that $d_\mathcal {T}(X,Z)\leq C_1$.

 \end{proposition}

\begin{proof}
By a theorem of Bers, there is a constant $c=c(S)$ such that every
$X\in \mathcal {M}(S)$ has a pants decomposition corresponding to a
maximal simplex $\sigma$ such that every curve of $\sigma$ has
length at most $c$ on $X$. With respect to these pants curves, each
of the twist coordinates is bounded by $2\pi c$, modulo the action
of Dehn twist about the curves in $\sigma$.

Now for a given $X$, we choose a point of $\Psi(\mathcal{V}(S))$
whose corresponding simplex has the topological type of $\sigma$.
For each curves $\alpha$ in $\sigma$ whose length is at most
$\epsilon_0$, we choose the corresponding Fenchel-Nielsen coordinate
of a point in $\Psi(\mathcal{V}(S))$ to be $l_X(\alpha)$. For each
curves $\beta$ in $\sigma$ whose length is between $\epsilon_0$ and
$c$, we choose the corresponding Fenchel-Nielsen coordinates of a
point of $\Psi(\mathcal{V}(S))$ to be $\epsilon_0$. In this way we
have chosen all the coordinates which determine a point $Z$ in
$\Psi(\mathcal{V}(S))$.

Since $X$ and $Z$ have bounded ratios of hyperbolic lengths and
bounded differences in twist coordinates, both $X$ and $Z$ either
project to a given thick part of moduli space or to the thin part.
In the first case, this implies the ratios of extremal lengths are
proportional to the ratios of hyperbolic lengths, and so one can use
Kerckhoff's distance formula to show that $d_\mathcal {T}(X,Z)$ are
bounded by some constant. In the second case one can use Minsky's
product theorem, noting that the constant $C$ in his formula is
universal, depending only on the topological type of $S$.
 \end{proof}

 \section{The proof of Theorem $1.3$}\label{s2}
Since $d_\mathcal {L}(X_1, X_2)\leq d_\mathcal {T}(X_1,X_2)$, by
Proposition $2.2$, we have shown that $\Psi
:(\mathcal{V}(S),d_{\mathcal{V}})\rightarrow
(\mathcal{M}(S),d_\mathcal {L})$ is almost onto . By Theorem $1.2$,
there is a constant D such that
$$d_\mathcal {L}(\Psi(x),\Psi(y))\leq d_\mathcal {T}(\Psi(x),\Psi(y))\leq d_\mathcal{V}(x,y) + D.$$
For the proof of Theorem $1.3$, it remains to show the opposite
inequality
\begin{equation}\label{equ:main}
d_\mathcal{V}(x,y)\leq d_\mathcal {L}(\Psi(x),\Psi(y)) + D'.
\end{equation}

Let $\sigma$ be a maximal simplex and  $P$ be the quotient map from $\mathcal{C}(S)$
to $\mathcal{C}(S)/\mathrm{Mod}(S)$. Denote the cone over $P(\sigma)$ by $\Delta(P(\sigma))$.

\begin{proposition}\label{pro:lt}
For any $Z_1, Z_2\in \Psi(\Delta(P(\sigma)))$, there is a constant
$C$, such that
$$d_\mathcal {T}(Z_1,Z_2)\leq d_\mathcal {L}(Z_1,Z_2) + C.$$
 \end{proposition}
\begin{proof}
Since $Z_i\in \Omega_\sigma(\epsilon_0)$ and the `` twist
coordinates '' of $Z_i$ are all vanishing, by Lemma $2.1$, we have
$$| d_\mathcal {T}(Z_1,Z_2) - \frac{1}{2}\log\sup_{\alpha_1,\cdots,\alpha_d}\frac{l_{Z_2}(\alpha_i)}{l_{Z_1}(\alpha_i)} |\leq C,$$
By the the length-spectra metric equation $(\ref{equ:ls})$, we have
\begin{eqnarray*}
d_\mathcal {T}(Z_1,Z_2)&\leq &\frac{1}{2}\log\sup_{\alpha_1,\cdots,\alpha_d}\frac{l_{Z_2}(\alpha_i)}{l_{Z_1}(\alpha_i)} + C\\
&\leq &d_\mathcal {L}(Z_1,Z_2) + C.
\end{eqnarray*}
 \end{proof}

We need the following technical lemma.
\begin{lemma} \label{lem:four}
There is a constant $C' = C'(S)$ such that any two points of
$\Psi(\mathcal{V}(S))$ can be joined by a $(1,C')$ quasi-geodesic in
the metric $d_\mathcal {L}$ that stays in $\Psi(\mathcal{V}(S))$ and
enters each simplex of $\Psi(\mathcal{V}(S))$ at most once.
 \end{lemma}

\begin{proof}
Denote the length of a path $\eta$ in $\mathcal{M}(S)$ with respect
to the Teichm\"{u}ller metric and the length-spectra metric by
$\|\eta\|_\mathcal {T}$ and $\|\eta\|_\mathcal {L}$.

As shown by Lemma 8 of Farb and Masur \cite{FM}, there is a constant
$C''$ such that if $\Psi(x), \Psi(y)$ lie in the same simplex
$\Psi(\Delta(P(\sigma)))$ of $\mathcal{M}(S)$, then there is a
$(1,C'')$ quasi-geodesic $\rho(x,y)$ in the metric $d_\mathcal {T}$
joining $\Psi(x)$ and $\Psi(y)$ that stays in
$\Psi(\Delta(P(\sigma)))$. By Proposition $3.1$,
\begin{eqnarray*}
\|\rho(x,y)\|_\mathcal {L} &\leq &
\|\rho(x,y)\|_\mathcal {T} \\
&\leq & d_\mathcal {T}(\Psi(x),\Psi(y))+ C''\\
&\leq & d_\mathcal {L}(\Psi(x),\Psi(y)) + C'' + C.\\
\end{eqnarray*}
 As a result, $\rho(x,y)$ is a
$(1,C''+C)$ quasi-geodesic in the metric $d_\mathcal {L}$.

Now suppose that $\Psi(x)\in \Psi(\Delta(P(\sigma_1)))$ and $\Psi(y)\in
\Psi(\Delta(P(\sigma_2)))$. If $\Psi(y)\in \Psi(\Delta(P(\sigma_1)))$ then
we are done by the argument above. Thus we can assume that
$\Psi(y)\not\in \Psi(\Delta(\sigma_1))$. Suppose that $\rho$ is a
geodesic about the length-spectrum metric from $\Psi(x)$ to
$\Psi(y)$ (the existence of geodesics about the length-spectrum
metric was shown by Thurston \cite{Thurston}, though maybe not
unique). Suppose $\rho$ leaves $\Psi(\Delta(P(\sigma_1)))$ and returns
to it for a last time at some $\Psi(z)\in
\Psi(\Delta(P(\sigma_1)))\bigcap \Psi(\Delta(P(\sigma_3)))$ for some
simplex $\Psi(\Delta(P(\sigma_3)))$. Then  we can replace $\rho$ by a
quasi-geodesic $\rho'$ that stays in $\Psi(\Delta(P(\sigma_1)))$ from $\Psi(x)$
to $\Psi(z)$ and then follows $\rho$ from $\Psi(z)$ to $\Psi(y)$
never returning to $ \Psi(\Delta(P(\sigma_1)))$. Note that the length
$\|\rho'\|_\mathcal{L}$ is less than $\|\rho\|_\mathcal{L}+C+C''$.

Continue with the above method, we now find the last
point $\Psi(w)$ that lies in $ \Psi(\Delta(P(\sigma_3)))$ and replace a
segment of $\rho'$ with the one that stays in $
\Psi(\Delta(P(\sigma_3)))$ and never returns again to $
\Psi(\Delta((\sigma_3)))$. Then we get a quasi-geodesic with length-spectra
length less than $\|\rho\|_\mathcal{L}+2(C+C'')$.

Since there are only a finite number of
maximal simplices in $\mathcal{C}(S)/\mathrm{Mod}(S)$, we can repeat the above operation
in a uniformly finite step. Then we get a path with length-spectra length larger
$d_\mathcal{L}(\Psi(x),\Psi(y))$ by an additive constant, which is a $(1,C')$ quasi-geodesic.
This prove the lemma.

 \end{proof}

We now continue with the final step in the proof of inequality
$(\ref{equ:main})$. Let $\rho$ be the $(1,C')$-quasi-geodesic in the
length-spectra metric as in Lemma $3.2$. Suppose that $\rho =
\bigcup_{i = 1,\cdots,n}{\rho_i}$ such that each $\rho_i\subset
\Psi(\Delta(P(\sigma_i)))$. By Theorem \ref{th:FM}, we have
\begin{eqnarray*}
d_{\mathcal{V}(S)}(x,y)&\leq &
d_\mathcal {T}(\Psi(x),\Psi(y)) +D \\
&\leq& \|\rho\|_\mathcal{T}+D \\
&\leq& \sum_{i=1,\cdots,n}\|\rho_i\|_
\mathcal {\mathcal {T}} +D\\
\end{eqnarray*}
It follows from the proof of Lemma \ref{lem:four} that each $\rho_i$ is also
a Teichm\"uller $(1,C'')$ quasi-geodesic. As a result, we know that the length
$\|\rho_i\|_\mathcal{T}$ is less than the Teichm\"uller distance between the two
endpoints of $\rho_i$. By Proposition \ref{pro:lt}, we have
$$\|\rho_i\|_\mathcal{T} \leq \|\rho_i\|_\mathcal{L}+C+C''.$$
Therefore,
\begin{eqnarray*}
d_{\mathcal{V}(S)}(x,y)&\leq &
 \sum_{i=1,\cdots,n}\|\rho_i\|_
\mathcal {\mathcal {L}} +D+n(C+C'')\\
&=& \|\rho\|_\mathcal{L}+D+n(C+C'') \\
&\leq& d_\mathcal {L}(\Psi(x),\Psi(y))+ C'+D+n(C+C'').
\end{eqnarray*}
As a result, we have proved the inequality $(\ref{equ:main})$.
\begin{remark}
Combining the known results, we can see that the Teichm\"{u}ller
metric, the length-spectra metric, Bergman metric, Carath\'{e}odory
mtric, McMullen metric, Ricci metric and perturbed Ricci metric (see
Liu, Sun and Yau \cite{LSY}) are quasi-isometric to each other on
moduli space. Note that Leuzinger [6], Farb and Weinberger [3]
proved that, while $\mathcal{M}(S)$ admits a metric of positive
scalar curvature for most $S$(when $g > 2$), it admits no metric of
positive scalar curvature with the same quasi-isometry type as the
Teichmuller metric on $\mathcal{M}(S)$.
\end{remark}

There is a natural question that whether the length-spectra metric and the Teichm\"uller metric
are bi-Lipschitz on moduli space. That is,
is there  a constant $K = K(S)$ such that
$d_\mathcal {L}(X_1,X_2)\leq d_\mathcal {T}(X_1,X_2)\leq K d_\mathcal {L}(X_1,X_2)$
for any $X_1, X_2\in \mathcal{M}(S)$ ?

Here the  left side inequality is trivial. If $d_\mathcal
{L}(X_1,X_2)$ is sufficiently large, the additive constant in
Theorem 1.4 can be absorbed into multiplicative constant to conclude
that $d_\mathcal {T} \leq K d_\mathcal {L}(X_1,X_2)$. Thus this
problem concerns mainly the local comparison. In a preprint \cite{ALPS}, we have
shown that such a constant $K$ depends on the injective radius of $X_1,X_2$.
There exists sequence $\{X_n, Y_n\}$ in the thin part of the Teichm\"uller space (or moduli space)
such that the ratio $d_{\mathrm{Teich}}(X_n, Y_n)/d_{\mathrm{ls}}(X_n, Y_n)$ tends to infinity.
As a result, the two metrics are not bi-Lipschitz in general.

\end{document}